# The probability of Riemann's hypothesis being true is equal to 1

Yuyang Zhu[1]

**Abstract** Let $P$ be the set of all prime numbers, $q_1, q_2, \cdots, q_m \in P$, $p_k$ be the $k$-th $(k = 1, 2, \cdots m)$ element of $P$ in ascending order of size, $\alpha_1, \alpha_2, \cdots, \alpha_m$ be positive integers, and $\beta_1, \beta_2, \cdots, \beta_m$ is a permutation of $\alpha_1, \alpha_2, \cdots, \alpha_m$ with $\beta_1 \geq \beta_2 \geq \cdots \geq \beta_m$. The following results are given in this paper:

(i) The following inequality is true

$$e^\gamma \log\log \prod_{k=1}^{m} q_k^{\alpha_k} - \prod_{k=1}^{m} \frac{q_k - \frac{1}{q_k^{\alpha_k}}}{q_k - 1} \geq e^\gamma \log\log \prod_{k=1}^{m} p_k^{\beta_k} - \prod_{k=1}^{m} \frac{p_k - \frac{1}{p_k^{\beta_k}}}{p_k - 1}.$$

(ii) If $n = \prod_{k=1}^{m} p_k^{\beta_k} = \left(\prod_{k=1}^{m} p_k\right)^{1+\varepsilon_m(n)}$, $\lim_{m \to \infty} \varepsilon_m(n) > 0$ or $\lim_{m \to \infty} \varepsilon_m(n) = +\infty$, then

$$\lim_{m \to \infty}(e^\gamma n \log\log n - \sigma(n)) > 0.$$

Where $\{\beta_k\}$ is a sequence, $\beta_k \in \mathbb{N}$, $\beta_1 \geq \beta_2 \geq \cdots \geq \beta_{m-1} \geq \beta_m$.

(iii) The probability of Riemann's hypothesis being true is equal to 1.

(iv) If $\lim_{m \to \infty} \varepsilon_m(n) = 0$ and there exist positive constant integers $A, K_0$, $A \geq \alpha_{K_0} \geq \alpha_{K_0+1} \geq \cdots \geq \alpha_{m-1} \geq \alpha_m$, then

$$\lim_{m \to \infty}(e^\gamma n \log\log n - \sigma(n)) > 0.$$

If $\lim_{m \to \infty} \varepsilon_m(n) = 0$, then there exist positive integers $K_1, K_2, \cdots, K_s (s \leq m)$, such that

$\alpha_i \leq K_i (i = 1, 2, \cdots, s)$ and $\lim_{m \to \infty} \frac{K}{m} = 0$. Where $K = \min\{K_1, K_2, \cdots, K_s\}$, $\sigma(n) = \sum_{d|n} d$, and $\gamma$ is the Euler constant.

**Keywords** Riemann hypothesis; non-trivial zeros; Robin's inequality; prime number.

[1] Department of Math. And Phys., Hefei University, Hefei China 230601. E-mail: zhuyy@hfuu.edu.cn.

**AMS Subject Classification 2000**   11M26

## 1 Introduction and main results

The Riemann hypothesis **[1]** (RH) states that the non-real zeros of the Riemann zeta function $\varsigma(s) = \sum_{n=1}^{+\infty} n^{-s}$ are on the line at $\Re(s) = 1/2$. Robin **[2]** pointed out that, if the RH is true, then $\sigma(n) \leqslant e^{\gamma} n \log \log n$ hold for all $n \geqslant 5041$, and $\gamma$ is Euler's constant. Conversely, If the RH is false, then there exist constants $0 < \beta < \frac{1}{2}$ and $C > 0$, such that

$$\sigma(n) \geqslant e^{\gamma} n \log \log n + \frac{Cn \log \log n}{(\log n)^{\beta}}$$

for infinite number of natural numbers $n$.

In 2007, Marek Wójtowicz **[3]** proved that, there is a subset $W$ of $\mathbb{N}$ of asymptotic density 1 such that $\lim_{n \to \infty, n \in W} k(n) = 0$, where $\mathbb{N}$ stands for the set of positive integers, and $k(n) = \frac{\sigma(n)}{e^{\gamma} n \log \log n}, (n > 2)$.

Zhu **[4-5]** proved the following inequality:
$$\lim_{n \to \infty} \frac{1}{n}\left((\exp H_n) \log H_n - \sigma(n)\right) \geqslant 0.$$

Solé and Zhu **[6]** proved the following conclusion:
$$\liminf_{n \to \infty} \left( e^{\gamma} \log \log n - \frac{1}{n} \sigma(n) \right) = 0.$$

Further, let $D(n) = e^{\gamma} n \log \log n - \sigma(n)$, we have the following limits when n ranges over Colossally Abundant Numbers **[6]**.

If RH is false, then $\liminf_{n \to \infty} D(n) = -\infty$.

If RH holds, then either $\lim_{n \to \infty} D(n) = \infty$, or $\liminf_{n \to \infty} D(n)$ is finite and $\geqslant 0$.

In fact, according to this conclusion and the Robin criterion for RH, the necessary and sufficient condition for RH being true is that

$$\lim_{n \to \infty} D(n) = \lim_{n \to \infty} (e^{\gamma} n \log \log n - \sigma(n)) \geqslant 0.$$

In this paper, our purpose is to further study the asymptotic property of Robin's inequality, and prove that the probability of RH being true is equal to 1.

### 1.1 Notations

For the convenience of the statement, we firstly give some notations.

**Notation 1** Denote a sequence of prime numbers by $q_k (k=1,2,\cdots,m)$ with $q_1 < q_2 < \cdots < q_m$.

**Notation 2** Denote by $P$ the set of all prime numbers and let $p_k$ be the $k$-th $(k=1,2,\cdots m)$ element of $P$ in ascending order of size.

### 1.2 Main results

The main results of the present paper are as follows.

**Theorem 1** Let $\alpha_k$ be a positive integer, $k=1,2,\cdots,m$, and $\beta_1,\beta_2,\cdots,\beta_m$ is a permutation of $\alpha_1,\alpha_2,\cdots,\alpha_m$ with $\beta_1 \geq \beta_2 \geq \cdots \geq \beta_m$, then we have

$$e^\gamma \log\log \prod_{k=1}^m q_k^{\alpha_k} - \prod_{k=1}^m \frac{q_k - \frac{1}{q_k^{\alpha_k}}}{q_k - 1} \geq e^\gamma \log\log \prod_{k=1}^m p_k^{\beta_k} - \prod_{k=1}^m \frac{p_k - \frac{1}{p_k^{\beta_k}}}{p_k - 1}. \qquad (1.1)$$

**Note 1** For any sufficiently large natural number $n' = \prod_{k=1}^m q_k^{\alpha_k}$, according to Theorem 1, if we prove inequality

$$l(n') = \frac{1}{n'} D(n') = \frac{1}{n'}(e^\gamma n' \log\log n' - \sigma(n'))$$

$$= e^\gamma \log\log \prod_{k=1}^m q_k^{\alpha_k} - \prod_{k=1}^m \frac{q_k - \frac{1}{q_k^{\alpha_k}}}{q_k - 1} > 0,$$

we only need to prove inequality

$$l(n) = \frac{1}{n} D(n) = \frac{1}{n}(e^\gamma n \log\log n - \sigma(n))$$

$$= e^\gamma \log\log \prod_{k=1}^m p_k^{\beta_k} - \prod_{k=1}^m \frac{p_k - \frac{1}{p_k^{\beta_k}}}{p_k - 1} > 0,$$

where $n = \prod_{k=1}^m p_k^{\beta_k}$. According to the Robin criterion for RH and the results in [6], we only need to study inequality

$$\lim_{n\to\infty}\frac{1}{n}D(n)=\lim_{n\to\infty}(e^{\gamma}\log\log\prod_{k=1}^{m}p_k^{\beta_k}-\prod_{k=1}^{m}\frac{p_k-\frac{1}{p_k^{\beta_k}}}{p_k-1})>0.$$

**Theorem 2** Let $\{\beta_k\}$ be a sequence, $\beta_k\in\mathbb{N}$, $\beta_1\geqslant\beta_2\geqslant\cdots\geqslant\beta_{m-1}\geqslant\beta_m$,

$n=\prod_{k=1}^{m}p_k^{\beta_k}=\left(\prod_{k=1}^{m}p_k\right)^{1+\varepsilon_m(n)}$, $P$ be the set of all prime numbers, $p_k$ be the k-th

$(k=1,2,\cdots m)$ element of $P$ in ascending order of size. If $\lim_{n\to\infty}\varepsilon_m(n)>0$ or

$\lim_{n\to\infty}\varepsilon_m(n)=+\infty$, then

$$\lim_{n\to\infty}(e^{\gamma}n\log\log n-\sigma(n))=+\infty. \qquad (1.2)$$

According to Theorem 1-2, we have

**Inference 1** There exists a positive constant $A$, such that for all natural numbers

$n=\prod_{k=1}^{m}q_k^{\alpha_k}$ with $2\leqslant\alpha_k\in\mathbb{N}$, $k=1,2,\cdots,m$, if $n\geqslant A$, then $D(n)>0$.

**Theorem 3** The probability of RH being true is equal to 1.

**Theorem 4** Let $\{\alpha_k\}$ be a sequence, $\alpha_k\in\mathbb{N}$, $\alpha_1\geqslant\alpha_2\geqslant\cdots\geqslant\alpha_{m-1}\geqslant\alpha_m$,

$n=\prod_{k=1}^{m}p_k^{\alpha_k}=\left(\prod_{k=1}^{m}p_k\right)^{1+\varepsilon_m(n)}$, $P$ be the set of all prime numbers, and $p_k$ be the k-th

$(k=1,2,\cdots m)$ element of $P$ in ascending order of size.

(1) If $\lim_{n\to\infty}\varepsilon_m(n)=0$ and there exist positive constant integers $A,K_0$, such

that $A\geqslant\alpha_{K_0}\geqslant\alpha_{K_0+1}\geqslant\cdots\geqslant\alpha_{m-1}\geqslant\alpha_m$, then

$$\lim_{n\to\infty}(e^{\gamma}n\log\log n-\sigma(n))>0.$$

(2) If $\lim_{n\to\infty}\varepsilon_m(n)=0$, then there exist positive integers $K_1,K_2,\cdots,K_s(s\leq m)$,

such that $\alpha_i\leq K_i(i=1,2,\cdots,s)$ and $\lim_{n\to\infty}\frac{K}{m}=0$, where $K=\min\{K_1,K_2,\cdots,K_s\}$.

According to Theorem 1 and (1) in the Theorem 4, we have

**Inference 2** There exists a positive constant $A$, for all natural numbers

$n=q_1q_2\cdots q_m$, if $n\geqslant A$, then $e^{\gamma}n\log\log n-\sigma(n)>0$.

**Conjecture** Let $\{\alpha_k\}$ be a sequence, $\alpha_k\in\mathbb{N}$, $\alpha_1\geqslant\alpha_2\geqslant\cdots\geqslant\alpha_{m-1}\geqslant\alpha_m$,

$$n = \prod_{k=1}^{m} p_k^{\alpha_k} = \left(\prod_{k=1}^{m} p_k\right)^{1+\varepsilon_m(n)},$$ $P$ be the set of all prime numbers, and $p_k$ be the k-th $(k=1,2,\cdots m)$ element of $P$ in ascending order of size.

If $\lim_{n\to\infty}\varepsilon_m(n) = 0$, then $\lim_{n\to\infty}(e^\gamma n\log\log n - \sigma(n)) > 0$.

If this conjecture is proved, then RH is also proved to be true.

## 1.3 Organization of article

We first give several important lemmas in Section 2. The proofs of main results are presented in Section 3. The proofs of lemmas are in Section 4.

## 2 Lemmas

In the following, we will give six lemmas as preliminaries for the proofs of Theorem 1~4.

For the convenience of statement, we appoint that the natural numbers are nonzero and $\mathbb{N}$ denotes the set of all natural numbers.

**Lemma 2.1** Let

$$f(q_1, q_2, \cdots, q_m) = e^\gamma \log\log \prod_{h=1}^{m} q_h^{\alpha_h} - \prod_{h=1}^{m} \frac{q_h - \frac{1}{q_h^{\alpha_h}}}{q_h - 1},$$

with $\alpha_k \in \mathbb{N}$, $k = 1, 2, \cdots, m$. Then the minimum of $f(q_1, q_2, \cdots, q_m)$ is

$$\min f(q_1, q_2, \cdots, q_m) = f(p_1, p_2, \cdots, p_m) = f(2, 3, 5, \cdots, p_m).$$

**Lemma 2.2([7, 8])** Let $\theta(p_m) = \sum_{k=1}^{m} \log p_k$, if $p_m > 10544111$, then

$$0.998684 p_m < \theta(p_m) < 1.001102 p_m. \tag{2.1}$$

$$p_m - \frac{0.0066788 p_m}{\log p_m} < \theta(p_m) < p_m + \frac{0.0066788 p_m}{\log p_m}. \tag{2.2}$$

Let $n = \prod_{k=1}^{m} p_k^{\alpha_k}$, where $\alpha_k \geq 1$, $\alpha_k \in \mathbb{N}$. $\varepsilon_m(n) = \left(\log \prod_{k=1}^{m} p_k^{\alpha_k}\right)\left(\log \prod_{k=1}^{m} p_k\right)^{-1} - 1$, then $\varepsilon_m(n) \geq 0$, and

$$\prod_{k=1}^{m} p_k^{\alpha_k} = \left(\prod_{k=1}^{m} p_k\right)^{1+\varepsilon_m(n)}.$$

We have the following Lemma 2.5:

**Lemma 2.3** Let $\{\alpha_k\}$ be a sequence, $\alpha_k \in \mathbb{N}$, and $\alpha_1 \geqslant \alpha_2 \geqslant \cdots \geqslant \alpha_{m-1} \geqslant \alpha_m$,

$$\prod_{k=1}^{m} p_k^{\alpha_k} = \left(\prod_{k=1}^{m} p_k\right)^{1+\varepsilon_m(n)}.$$

If $\lim\limits_{m\to\infty} \varepsilon_m(n) \neq 0$, then there exists constant $a > 0$, such that $\lim\limits_{m\to\infty} \varepsilon_m(n) > 2a > 0$.

**Lemma 2.4** The following inequality holds:

$$\lim_{m\to+\infty}\left(e^\gamma \log\log \prod_{k=1}^{m} p_k - \prod_{k=1}^{m} \frac{p_k - \frac{1}{p_k^{\alpha_k}}}{p_k - 1}\right) \geqslant 0. \qquad (2.3)$$

**Lemma 2.5([9])** $p_k = k\{\log k + \log\log k - 1 + O(\frac{\log\log k}{\log k})\}$.

**Lemma 2.6[6]** We have the following limits when n ranges over Colossally Abundant Numbers.

If RH is false, then $\liminf\limits_{n\to\infty} D(n) = -\infty$.

If RH holds, then either $\lim\limits_{n\to\infty} D(n) = \infty$, or $\liminf\limits_{n\to\infty} D(n)$ is finite and $\geqslant 0$.

# 3 Proofs of Theorems

## 3.1 Proof of Theorem 1

To prove Theorem 1, the following propositions are required.

**Proposition 3.1** If $1 \leqslant \alpha \leqslant \beta$, $1 \leqslant a \leqslant b$, then

$$\left(1 - \frac{1}{a^\alpha}\right)\left(1 - \frac{1}{b^\beta}\right) \leqslant \left(1 - \frac{1}{a^\beta}\right)\left(1 - \frac{1}{b^\alpha}\right). \qquad (3.1)$$

**Proof** Since $1 \leqslant \alpha \leqslant \beta$, $1 \leqslant a \leqslant b$, (3.1) is equivalent to

$$\left(1 - \frac{1}{a^\alpha}\right)\left(1 - \frac{1}{a^\beta}\right)^{-1} \leqslant \left(1 - \frac{1}{b^\alpha}\right)\left(1 - \frac{1}{b^\beta}\right)^{-1}. \qquad (3.2)$$

Let

$$v(x) = \left(1 - \frac{1}{x^\alpha}\right)\left(1 - \frac{1}{x^\beta}\right)^{-1} \quad (x > 1).$$

We will prove $v(x)$ is a monotonically increasing function. Since

$$v'(x) = \frac{\alpha}{x^{\alpha+1}}\left(1-\frac{1}{x^\beta}\right)^{-1} - \frac{\beta}{x^{\beta+1}}\left(1-\frac{1}{x^\alpha}\right)\left(1-\frac{1}{x^\beta}\right)^{-2}$$

$$= \frac{1}{x^{\alpha+\beta+1}}\left(1-\frac{1}{x^\beta}\right)^{-2}(\alpha x^\beta - \beta x^\alpha + \beta - \alpha),$$

and

$$x>1,\ 1 \leqslant \alpha \leqslant \beta,\ \frac{1}{x^{\alpha+\beta+1}}\left(1-\frac{1}{x^\beta}\right)^{-2} > 0,$$

the proof of $v'(x) \geqslant 0$ is equivalent to proving that $\alpha x^\beta - \beta x^\alpha + \beta - \alpha$ is non-negative.

Let $w(x) = \alpha x^\beta - \beta x^\alpha + \beta - \alpha$. Then

$$w'(x) = \alpha\beta(x^{\beta-1} - x^{\alpha-1}).$$

Since

$$x \geqslant 1,\text{ and } 1 \leqslant \alpha \leqslant \beta,\ w'(x) = \alpha\beta(x^{\beta-1} - x^{\alpha-1}) \geqslant 0,$$

which shows that $w(x)$ is a monotonically increasing function in $[1,+\infty)$. In addition, $w(x)$ is continuous at $x=1$. Therefore, we can conclude that for $x \geqslant 1$,

$$w(x) \geqslant w(1) = 0.$$

Hence $\alpha x^\beta - \beta x^\alpha + \beta - \alpha \geqslant 0$. Furthermore, $v'(x) \geqslant 0$, $v(x)$ is a monotonically increasing function. (3.1) and (3.2) are proved. □

According to Proposition 3.1, the following inequality is obvious,

$$\left(1-\frac{1}{a^{\alpha+1}}\right)\left(1-\frac{1}{b^{\beta+1}}\right) \leqslant \left(1-\frac{1}{a^{\beta+1}}\right)\left(1-\frac{1}{b^{\alpha+1}}\right), \qquad (3.3)$$

with $1 \leqslant \alpha \leqslant \beta$ and $1 \leqslant a \leqslant b$.

According to (3.3),

$$ab\left(1-\frac{1}{a^{\alpha+1}}\right)\left(1-\frac{1}{b^{\beta+1}}\right) \leqslant ab\left(1-\frac{1}{a^{\beta+1}}\right)\left(1-\frac{1}{b^{\alpha+1}}\right),$$

namely

$$\left(a-\frac{1}{a^\alpha}\right)\left(b-\frac{1}{b^\beta}\right) \leqslant \left(a-\frac{1}{a^\beta}\right)\left(b-\frac{1}{b^\alpha}\right). \qquad (3.4)$$

Consequently, we have

**Proposition 3.2** Let $\alpha_i \geqslant 1\ (i=1,2,\cdots,m)$, $1 < q_1 \leqslant q_2 \leqslant \cdots \leqslant q_m$, $\beta_1, \beta_2, \cdots, \beta_m$ be a permutation of $\alpha_1, \alpha_2, \cdots, \alpha_m$, and $\beta_1 \geqslant \beta_2 \geqslant \cdots \geqslant \beta_m$. Then

$$\prod_{i=1}^{m}\left(q_i - \frac{1}{q_i^{\alpha_i}}\right) \leqslant \prod_{i=1}^{m}\left(q_i - \frac{1}{q_i^{\beta_i}}\right). \qquad (3.5)$$

**Proof** For $m=2$, according to (3.4), this proposition is true. Assume the proposition holds for $m=k\ (k \geqslant 2)$. Then for $m=k+1$, let $\beta_1 = \max\{\alpha_1, \alpha_2, \cdots, \alpha_{k+1}\} = \alpha_r\ (1 \leqslant r \leqslant k+1)$. We will prove the proposition 3.2 in two cases: ① $\beta_1 \neq \alpha_1$, ② $\beta_1 = \alpha_1$.

① If $\beta_1 \neq \alpha_1$, then $\beta_1 = \alpha_r$, $2 \leqslant r \leqslant k+1$. According to (3.4),

$$\prod_{i=1}^{k+1}\left(q_i - \frac{1}{q_i^{\alpha_i}}\right) = \left(q_1 - \frac{1}{q_1^{\alpha_1}}\right)\left(q_r - \frac{1}{q_r^{\alpha_r}}\right)\prod_{2 \leq i \leq k+1, i \neq r}\left(q_i - \frac{1}{q_i^{\alpha_i}}\right)$$

$$\leqslant \left(q_1 - \frac{1}{q_1^{\alpha_r}}\right)\left(q_r - \frac{1}{q_r^{\alpha_1}}\right)\prod_{2 \leq i \leq k+1, i \neq r}\left(q_i - \frac{1}{q_i^{\alpha_i}}\right)$$

$$= \left(q_1 - \frac{1}{q_1^{\beta_1}}\right)\left[\left(q_r - \frac{1}{q_r^{\alpha_1}}\right)\prod_{2 \leq i \leq k+1, i \neq r}\left(q_i - \frac{1}{q_i^{\alpha_i}}\right)\right]. \qquad (3.6)$$

According to the assumption of $m=k$,

$$\left(q_r - \frac{1}{q_r^{\alpha_1}}\right)\prod_{2 \leq i \leq k+1, i \neq r}\left(q_i - \frac{1}{q_i^{\alpha_i}}\right) \leqslant \prod_{i=2}^{k+1}\left(q_i - \frac{1}{q_i^{\beta_i}}\right). \qquad (3.7)$$

By (3.6) and (3.7),

$$\prod_{i=1}^{k+1}\left(q_i - \frac{1}{q_i^{\alpha_i}}\right) \leqslant \prod_{i=1}^{k+1}\left(q_i - \frac{1}{q_i^{\beta_i}}\right).$$

② If $\beta_1 = \alpha_1$, according to the assumption for $m=k$, i.e.,

$$\prod_{i=2}^{k+1}\left(q_i - \frac{1}{q_i^{\alpha_i}}\right) \leqslant \prod_{i=2}^{k+1}\left(q_i - \frac{1}{q_i^{\beta_i}}\right). \qquad (3.8)$$

By (3.8),

$$\prod_{i=1}^{k+1}\left(q_i - \frac{1}{q_i^{\alpha_i}}\right) = \left(q_1 - \frac{1}{q_1^{\alpha_1}}\right)\prod_{i=2}^{k+1}\left(q_i - \frac{1}{q_i^{\alpha_i}}\right) = \left(q_1 - \frac{1}{q_1^{\beta_1}}\right)\prod_{i=2}^{k+1}\left(q_i - \frac{1}{q_i^{\alpha_i}}\right)$$

$$\leqslant \left(q_1 - \frac{1}{q_1^{\beta_1}}\right)\prod_{i=2}^{k+1}\left(q_i - \frac{1}{q_i^{\beta_i}}\right) = \prod_{i=1}^{k+1}\left(q_i - \frac{1}{q_i^{\beta_i}}\right).$$

To conclude, the proposition is true for $m=k+1$. It is proved by the principle of mathematical induction. □

**Proposition 3.3** If $1 \leqslant \alpha \leqslant \beta$, $1 \leqslant a \leqslant b$, then

$$a^\alpha b^\beta \geqslant a^\beta b^\alpha.$$

**Proof** From $a^{\beta-\alpha} \leqslant b^{\beta-\alpha}$, $a^\alpha b^\beta \geqslant a^\beta b^\alpha$. □

From Proposition 3.3 and the principle of mathematical induction, we can prove the following conclusion:

**Proposition 3.4** Let $\alpha_i \geqslant 1$ $(i=1,2,\cdots,m)$, $1 < q_1 \leqslant q_2 \leqslant \cdots \leqslant q_m$, $\beta_1, \beta_2, \cdots, \beta_m$ be a permutation of $\alpha_1, \alpha_2, \cdots, \alpha_m$, and $\beta_1 \geqslant \beta_2 \geqslant \cdots \geqslant \beta_m$. Then

$$\prod_{k=1}^m q_k^{\alpha_k} \geqslant \prod_{k=1}^m q_k^{\beta_k}. \tag{3.9}$$

**Proof** $m=1$ is a trivial case. For $m=2$, the proposition is true according to Proposition 3.3. Assume the proposition is true for $m=t \geqslant 2$. Then for $m=t+1$, let $\beta_{t+1} = \min_{1 \leqslant k \leqslant t+1}\{\alpha_k\} = \alpha_j$ $(1 \leqslant j \leqslant t+1)$. If $j \neq t+1$, according to Proposition 3.3,

$$q_j^{\alpha_j} q_{t+1}^{\alpha_{t+1}} \geqslant q_j^{\alpha_{t+1}} q_{t+1}^{\alpha_j} = q_j^{\alpha_{t+1}} q_{t+1}^{\beta_{t+1}},$$

$$\prod_{k=1}^{t+1} q_k^{\alpha_k} = q_j^{\alpha_j} q_{t+1}^{\alpha_{t+1}} \prod_{k \neq j, t+1} q_k^{\alpha_k}$$

$$\geqslant q_j^{\alpha_{t+1}} q_{t+1}^{\beta_{t+1}} \prod_{k \neq j, t+1} q_k^{\alpha_k} = q_{t+1}^{\beta_{t+1}}\left(q_j^{\alpha_{t+1}} \prod_{k \neq j, t+1} q_k^{\alpha_k}\right).$$

Using the induction hypothesis when $m=t$, i.e.,

$$q_j^{\alpha_{t+1}} \prod_{k \neq j, t+1} q_k^{\alpha_k} \geqslant \prod_{k=1}^t q_k^{\beta_k},$$

(3.9) holds. Thus

$$\prod_{k=1}^{t+1} q_k^{\alpha_k} \geqslant q_{t+1}^{\beta_{t+1}}\left(q_j^{\alpha_{t+1}} \prod_{k \neq j, t+1} q_k^{\alpha_k}\right) \geqslant q_{t+1}^{\beta_{t+1}} \prod_{k=1}^t q_k^{\beta_k} = \prod_{k=1}^{t+1} q_k^{\beta_k}.$$

Hence, (3.9) holds for $j \neq t+1$.

If $j=t+1$, then $\beta_{t+1} = \min_{1 \leqslant k \leqslant t+1}\{\alpha_k\} = \alpha_{t+1}$, according to the induction hypothesis when $m=t$,
$$\prod_{k=1}^{t} q_k^{\alpha_k} \geqslant \prod_{k=1}^{t} q_k^{\beta_k},$$
then
$$\prod_{k=1}^{t+1} q_k^{\alpha_k} = q_{t+1}^{\alpha_{t+1}} \prod_{k=1}^{t} q_k^{\alpha_k} \geqslant q_{t+1}^{\alpha_{t+1}} \prod_{k=1}^{t} q_k^{\beta_k}$$
$$= q_{t+1}^{\beta_{t+1}} \prod_{k=1}^{t} q_k^{\beta_k} = \prod_{k=1}^{t+1} q_k^{\beta_k}.$$

Since $\beta_{t+1} = \min_{1 \leq k \leq t+1}\{\alpha_k\} = \alpha_{t+1}$, $\beta_1 \geqslant \beta_2 \geqslant \cdots \geqslant \beta_t \geqslant \beta_{t+1}$, and $\beta_1, \beta_2, \cdots, \beta_t, \beta_{t+1}$ is a permutation of $\alpha_1, \alpha_2, \cdots, \alpha_t, \alpha_{t+1}$. Hence, (3.9) holds for $j=t+1$. The proposition is true for $m=t+1$. From the principle of mathematical induction, the proposition is true. □

Now, the proof of theorem 1 can be given as follows:

**Proof of Theorem 1** According to (3.9), we have
$$e^{\gamma} \log\log \prod_{k=1}^{m} q_k^{\alpha_k} \geqslant e^{\gamma} \log\log \prod_{k=1}^{m} q_k^{\beta_k}. \tag{3.10}$$

According to (3.5),
$$\prod_{k=1}^{m} \frac{q_k - \frac{1}{q_k^{\alpha_k}}}{q_k - 1} = \left(\prod_{k=1}^{m} \frac{1}{q_k - 1}\right)\prod_{k=1}^{m}\left(q_k - \frac{1}{q_k^{\alpha_k}}\right)$$
$$\leqslant \left(\prod_{k=1}^{m} \frac{1}{q_k - 1}\right)\prod_{k=1}^{m}\left(q_k - \frac{1}{q_k^{\beta_k}}\right) = \prod_{k=1}^{m} \frac{q_k - \frac{1}{q_k^{\beta_k}}}{q_k - 1}. \tag{3.11}$$

By (3.10) and (3.11),
$$e^{\gamma} \log\log \prod_{k=1}^{m} q_k^{\alpha_k} - \prod_{k=1}^{m} \frac{q_k - \frac{1}{q_k^{\alpha_k}}}{q_k - 1} \geqslant e^{\gamma} \log\log \prod_{k=1}^{m} q_k^{\beta_k} - \prod_{k=1}^{m} \frac{q_k - \frac{1}{q_k^{\beta_k}}}{q_k - 1}. \tag{3.12}$$

According to Lemma 2.1,
$$e^{\gamma} \log\log \prod_{k=1}^{m} q_k^{\beta_k} - \prod_{k=1}^{m} \frac{q_k - \frac{1}{q_k^{\beta_k}}}{q_k - 1} \geqslant e^{\gamma} \log\log \prod_{k=1}^{m} p_k^{\beta_k} - \prod_{k=1}^{m} \frac{p_k - \frac{1}{p_k^{\beta_k}}}{p_k - 1}. \tag{3.13}$$

Thus (1.1) is true by (3.12) and (3.13). □

## 3.2 Proof of Theorem 2

**Proof of Theorem 2**

Since $n = \prod_{k=1}^{m} p_k^{\beta_k} = \left(\prod_{k=1}^{m} p_k\right)^{1+\varepsilon_m(n)}$, $\beta_1 \geq \beta_2 \geq \cdots \geq \beta_m \geq 1$, if $m$ is finite, then there exists $N_0 \in \mathbb{N}$, such that $m < N_0$. When $n \to +\infty$, then $\varepsilon_m(n) \to +\infty$ (from $\prod_{k=1}^{N_0} p_k$ is finite) and $e^\gamma \log\log n \to +\infty$. Since $N_0$ is a positive constant, $m$ is finite,

$$\prod_{k=1}^{m} \frac{q_k - \frac{1}{q_k^{\alpha_k}}}{q_k - 1} < 2^m < 2^{N_0}.$$

Hence,

$$\lim_{n \to \infty} \frac{1}{n}\left(e^\gamma n \log\log n - \sigma(n)\right)$$

$$= \lim_{n \to +\infty} \frac{1}{n}\left(e^\gamma n \log\log n - n\prod_{k=1}^{m} \frac{q_k - \frac{1}{q_k^{\alpha_k}}}{q_k - 1}\right)$$

$$\geq \lim_{n \to +\infty}\left(e^\gamma \log(\log n) - 2^{N_0}\right) = +\infty.$$

namely,

$$\lim_{n \to +\infty}(e^\gamma n \log\log n - \sigma(n)) = +\infty.$$

That is (1.2) is proved under the assumption of $m < N_0$.

Now we need to consider the $m \to \infty$ case. According to Lemma 2.3, there exists constant $a > 0$, such that $\lim_{m \to \infty} \varepsilon_m(n) \geq 2a > 0$.

In addition,

$$e^\gamma \log\log \prod_{k=1}^{m} p_k^{\beta_k} = e^\gamma \log\log\left(\prod_{k=1}^{m} p_k\right)^{1+\varepsilon_m(n)}$$

$$= e^\gamma \log(1+\varepsilon_m(n)) + e^\gamma \log\log \prod_{k=1}^{m} p_k,$$

namely

$$\frac{1}{n}\left(e^\gamma n \log\log n - \sigma(n)\right)$$

$$\geq = e^\gamma \log\log \prod_{k=1}^{m} p_k^{\beta_k} - \prod_{k=1}^{m} \frac{p_k - \frac{1}{p_k^{\beta_k}}}{p_k - 1}$$

$$= e^{\gamma}\log(1+\varepsilon_m(n)) + e^{\gamma}\log\log\left(\prod_{k=1}^{m} p_k\right) - \prod_{k=1}^{m}\frac{p_k - \frac{1}{p_k^{\beta_k}}}{p_k - 1}.$$

When $m \to +\infty$, according to Lemma 2.3 and 2.4,

$$\lim_{n\to\infty}\frac{1}{n}\left(e^{\gamma}n\log\log n - \sigma(n)\right)$$

$$= \lim_{m\to\infty}\left(e^{\gamma}\log\log\prod_{k=1}^{m} p_k^{\beta_k} - \prod_{k=1}^{m}\frac{p_k - \frac{1}{p_k^{\beta_k}}}{p_k - 1}\right)$$

$$= \lim_{m\to\infty}\left(e^{\gamma}\log(1+\varepsilon_m(n)) + e^{\gamma}\log\log\left(\prod_{k=1}^{m} p_k\right) - \prod_{k=1}^{m}\frac{p_k - \frac{1}{p_k^{\beta_k}}}{p_k - 1}\right)$$

$$= \lim_{m\to\infty} e^{\gamma}\log(1+\varepsilon_m(n)) + \lim_{m\to\infty}\left(e^{\gamma}\log\log\left(\prod_{k=1}^{m} p_k\right) - \prod_{k=1}^{m}\frac{p_k - \frac{1}{p_k^{\beta_k}}}{p_k - 1}\right)$$

$$\geq e^{\gamma}\log(1+2a) + 0 > e^{\gamma}\log(1+a) > 0.$$

That is

$$\lim_{n\to+\infty}(e^{\gamma}n\log\log n - \sigma(n)) = +\infty \geq \lim_{n\to+\infty} ne^{\gamma}\log(1+a) = +\infty.$$

Thus (1.2) holds.

In summary, when $\lim_{n\to\infty}\varepsilon_m(n) > 0$ or $\lim_{n\to\infty}\varepsilon_m(n) = +\infty$, then

$$\lim_{n\to+\infty}(e^{\gamma}n\log\log n - \sigma(n)) = +\infty. \square$$

**Proof of inference 1** since $\alpha_k \geq 2, k = 1, 2, \cdots, m$, $\varepsilon_m(n) \geq 1$, $\lim_{n\to\infty}\varepsilon_m(n) \geq 1$, according to Theorem 1 and Theorem 2, we have $\lim_{n\to+\infty}(e^{\gamma}n\log\log n - \sigma(n)) = +\infty$, namely there exists a positive constant $A$, such that for all natural number, if $n \geq A$, then $D(n) = e^{\gamma}n\log\log n - \sigma(n) > 0$. $\square$

### 3.3 Proof of Theorem 3

**Proof of Theorem 3** Let $n = \prod_{k=1}^{m} q_k^{\alpha_k}$ ($q_1, q_2, \cdots, q_m \in P, \alpha_1, \alpha_2, \cdots, \alpha_m \in \mathbb{N}$), by Theorem 1, we only need to consider $n = \prod_{k=1}^{m} p_k^{\beta_k} = \left(\prod_{k=1}^{m} p_k\right)^{1+\varepsilon_m(n)}$, $\beta_k \in \mathbb{N}$, $\beta_1 \geq \beta_2 \geq \cdots \geq \beta_m$, therefore $\varepsilon_m(n) = \left(\log\prod_{k=1}^{m} p_k^{\beta_k}\right)\left(\log\prod_{k=1}^{m} p_k\right)^{-1} - 1$ and $\varepsilon_m(n) \geq 0$. If

$\varepsilon_m(n) \in [0,1]$ then for $n_t = \prod_{k=1}^{m} p_k^{t+\beta_k}, (t \in \mathbb{N})$,

$$\varepsilon_m(n_t) = \left(\log \prod_{k=1}^{m} p_k^{t+\beta_k}\right)\left(\log \prod_{k=1}^{m} p_k\right)^{-1} - 1$$

$$= \left(\log \prod_{k=1}^{m} p_k^{t} + \log \prod_{k=1}^{m} p_k^{\beta_k}\right)\left(\log \prod_{k=1}^{m} p_k\right)^{-1} - 1$$

$$= t + \left(\log \prod_{k=1}^{m} p_k^{\beta_k}\right)\left(\log \prod_{k=1}^{m} p_k\right)^{-1} - 1 = t + \varepsilon_m(n).$$

$\varepsilon_m(n_t) = t + \varepsilon_m(n) \in [t, t+1]$. Similarly, if $\varepsilon_m(n) \in [t, t+1]$ $(t \in \mathbb{N})$, then for $n_{-t} = \prod_{k=1}^{m} p_k^{\beta_k - t}$, $\varepsilon_m(n_{-t}) = \varepsilon_m(n) - t \in [0,1]$, since $t$ is arbitrary, $\varepsilon_m(n)$ distributes evenly in every interval $[t, t+1](t = 0, 1, 2, \cdots)$. According to Theorem 2, if $\lim_{n \to \infty} \varepsilon_m(n) > 0$ or $\lim_{n \to \infty} \varepsilon_m(n) = +\infty$, then

$$\lim_{n \to \infty}\left(e^{\gamma} n \log \log n - \sigma(n)\right) = +\infty,$$

hence, if $\lim_{n \to \infty}\left(e^{\gamma} n \log \log n - \sigma(n)\right) \leq 0$, then $\lim_{n \to \infty} \varepsilon_m(n) = 0$. According to $\lim_{m \to \infty} \varepsilon_m(n) \geq 0$ and Lemma 2.6, the probability that RH being true is equal to 1.□

**Note 2**: An event with a probability equal to 1 is not necessarily an inevitable event. According to Theorem 2 and Lemma 2.6, when $\lim_{m \to \infty} \varepsilon_m(n) = 0$, (1.2) hold, then RH is true. Theorem 4 discusses the asymptotic property of Robin inequality when $\lim_{m \to \infty} \varepsilon_m(n) = 0$.

### 3.4 Proof of Theorem 4

(1) Let $R(x) = e^{\gamma} \log \log \prod_{k=1}^{m} p_k^{\alpha_k x} - \prod_{k=1}^{m} \dfrac{p_k - \dfrac{1}{p_k^{\alpha_k x}}}{p_k - 1}$, $x \in [1, 2]$

$$R'(x) = \frac{e^{\gamma}}{x} - \left[\prod_{k=1}^{m} \frac{p_k - \left(\frac{1}{p_k^{\alpha_k}}\right)^x}{p_k - 1}\right] \sum_{k=1}^{m} \frac{\alpha_k \log p_k}{p_k^{x\alpha_k+1} - 1}$$

$$< e^{\gamma} - \frac{\alpha_{K_0} \log p_{K_0}}{p_{N_0}^{x\alpha_{K_0}+1} - 1} \prod_{k=1}^{m} \frac{p_k - \left(\frac{1}{p_k^{\alpha_k}}\right)^x}{p_k - 1}$$

$$\leqslant e^{\gamma} - \frac{\log p_{K_0}}{p_{K_0}^{Ax+1} - 1} \prod_{k=1}^{m} \frac{p_k - \left(\frac{1}{p_k^{\alpha_k}}\right)^x}{p_k - 1}$$

$$\leqslant e^{\gamma} - \frac{\log p_{K_0}}{p_{K_0}^{Ax+1} - 1} \prod_{k=1}^{m} \frac{p_k - \frac{1}{p_k}}{p_k - 1} = e^{\gamma} - \frac{\log p_{K_0}}{p_{K_0}^{Ax+1} - 1} \prod_{k=1}^{m} \left(\frac{p_k + 1}{p_k}\right).$$

But $\lim_{m \to \infty} \prod_{k=1}^{m} \left(\frac{p_k + 1}{p_k}\right) = +\infty$, hence

$$\lim_{m \to \infty} \left[e^{\gamma} - \frac{\log p_{K_0}}{p_{K_0}^{Ax+1} - 1} \prod_{k=1}^{m} \left(\frac{p_k + 1}{p_k}\right)\right] = -\infty.$$

Namely $\lim_{m \to \infty} R'(x) = -\infty$, then there exists a positive integer $N$, when $m > N$,

we have $R(1) > R(2)$. Thus for $m > N$,

$$R(1) = e^{\gamma} \log \log \prod_{k=1}^{m} p_k^{\alpha_k} - \prod_{k=1}^{m} \frac{p_k - \frac{1}{p_k^{\alpha_k}}}{p_k - 1}$$

$$> R(2) = e^{\gamma} \log \log \left(\prod_{k=1}^{m} p_k^{\alpha_k}\right)^2 - \prod_{k=1}^{m} \frac{p_k - \left(\frac{1}{p_k^{\alpha_k}}\right)^2}{p_k - 1}$$

$$= e^{\gamma} \log \log \prod_{k=1}^{m} p_k^{2\alpha_k} - \prod_{k=1}^{m} \frac{p_k - \frac{1}{p_k^{2\alpha_k}}}{p_k - 1}.$$

Let $\omega_k = 2\alpha_k$, then $\omega_k \geqslant 2$ $(k = 1, 2, \cdots, m)$. For a natural number $\bar{n} = n^2 = \prod_{k=1}^{m} p_k^{\omega_k}$, since $\omega_k \geqslant 2$, $\varepsilon_m(\bar{n}) \geq 1 > 0$ by Theorem 2,

$$\lim_{\bar{n} \to \infty} R(2) = \lim_{\bar{n} \to \infty} \frac{1}{\bar{n}} \left(e^{\gamma} \bar{n} \log \log \bar{n} - \sigma(\bar{n})\right)$$

$$= \lim_{\overline{n}\to\infty}(e^{\gamma}\log\log\prod_{k=1}^{m}p_k^{2\alpha_k}-\prod_{k=1}^{m}\frac{p_k-\frac{1}{p_k^{2\alpha_k}}}{p_k-1})>0.$$

$$\lim_{n\to\infty}R(1)=\lim_{n\to\infty}\frac{1}{n}\left(e^{\gamma}n\log\log n-\sigma(n)\right)$$

$$=\lim_{n\to\infty}(e^{\gamma}\log\log\prod_{k=1}^{m}p_k^{\alpha_k}-\prod_{k=1}^{m}\frac{p_k-\frac{1}{p_k^{\alpha_k}}}{p_k-1})\geq\lim_{\overline{n}\to\infty}R(2).$$

Thus, $\lim_{n\to\infty}\frac{1}{n}\left(e^{\gamma}n\log\log n-\sigma(n)\right)>0$, namely $\lim_{n\to\infty}\left(e^{\gamma}n\log\log n-\sigma(n)\right)=+\infty$.

**(2)** According to $\lim_{m\to\infty}\varepsilon_m(n)=0$, $\alpha_1\geq\alpha_2\geq\cdots\geq\alpha_m\geq 1$, there exist positive integers $K_1,K_2,\cdots,K_s(s\leq m)$, such that $\alpha_{K_i}\leq K_i(i=1,2,\cdots,s)$ ( if not, then $\forall j\in\{1,2,\cdots,m\}$, $\alpha_j>j$, thus $\varepsilon_m(n)\geq 1$, $\lim_{m\to\infty}\varepsilon_m(n)\geq 1$, but $\lim_{m\to\infty}\varepsilon_m(n)=0$, which is contradictory ).

Let $K_0=\min\{K_1,K_2,\cdots,K_s\}$, then $K_0\geq 1$. When $K_0=1$, then $\lim_{m\to\infty}\frac{K_0}{m}=\lim_{m\to\infty}\frac{1}{m}=0$. If $K_0\geq 2$ we need to prove that $\lim_{m\to\infty}\frac{K_0}{m}=0$.

If $\lim_{m\to\infty}\frac{K_0}{m}=b>0$, then $K_0=bm+r_m$ (with $\lim_{m\to\infty}r_m=0$). Since

$$n=\prod_{k=1}^{m}p_k^{\alpha_k}=\left(\prod_{k=1}^{m}p_k\right)^{1+\varepsilon_m(n)}\geq(\prod_{k=1}^{K_0}p_k^{K_0})\prod_{k=K_0+1}^{m}p_k=(\prod_{k=1}^{K_0}p_k^{K_0-1})\prod_{k=1}^{m}p_k.$$

Namely

$$\left(\prod_{k=1}^{m}p_k\right)^{1+\varepsilon_m(n)}\geq(\prod_{k=1}^{K_0}p_k^{K_0-1})\prod_{k=1}^{m}p_k, \text{ thus } \left(\prod_{k=1}^{m}p_k\right)^{\varepsilon_m(n)}\geq\prod_{k=1}^{K_0}p_k^{K_0-1}.$$

Taking logarithm on both sides

$$\varepsilon_m(n)\theta(p_m)=\varepsilon_m(n)\sum_{k=1}^{m}\log p_k\geq(K_0-1)\sum_{k=1}^{K_0}\log p_k=(K_0-1)\theta(p_{K_0}).$$

By (2.2) in Lemma 2.2, when $p_{K_0}>10544111$,

$$\varepsilon_m(n)p_m(1+\frac{0.0066788}{\log p_m})>\varepsilon_m(n)\theta(p_m)\geq(K_0-1)\theta(p_{K_0})$$

$$>(bm+r_m-1)p_{bm+r_m}(1-\frac{0.0066788}{\log p_{bm+r_m}}). \qquad (3.14)$$

According to Lemma 2.5 and (3.14),

$$\varepsilon_m(n)\left[m\{\log m+\log\log m-1+O\left(\frac{\log\log m}{\log m}\right)\}\right](1+\frac{0.0066788}{\log p_m})>$$

$$(bm+r_m-1)(bm+r_m)\{\log(bm+r_m)+(\log\log(bm+r_m)-1)$$

$$+O(\frac{\log\log(bm+r_m)}{\log(bm+r_m)})\}(1-\frac{0.0066788}{\log p_{bm+r_m}}).$$

Dividing both sides by $m\log m$,

$$\varepsilon_m(n)\left[\{1+\frac{\log\log m-1}{\log m}+O\left(\frac{\log\log m}{(\log m)^2}\right)\}\right](1+\frac{0.0066788}{\log p_m})>$$

$$(b\frac{m}{\log m}+\frac{r_m-1}{\log m})(b+\frac{r_m}{m})\{\log(bm+r_m)+(\log\log(bm+r_m)-1)$$

$$+O\left(\frac{\log\log(bm+r_m)}{\log(bm+r_m)}\right)\}(1-\frac{0.0066788}{\log p_{bm+r_m}}).$$

But

$$\lim_{m\to\infty}\varepsilon_m(n)\left[\{1+\frac{\log\log m-1}{\log m}+O\left(\frac{\log\log m}{(\log m)^2}\right)\}\right](1+\frac{0.0066788}{\log p_m})$$

$$=\lim_{m\to\infty}\varepsilon_m(n)=0.$$

$$\lim_{m\to\infty}(b\frac{m}{\log m}+\frac{r_m-1}{\log m})(b+\frac{r_m}{m})\{\log(bm+r_m)+(\log\log(bm+r_m)-1)$$

$$+O(\frac{\log\log(bm+r_m)}{\log(bm+r_m)})\}(1-\frac{0.0066788}{\log p_{bm+r_m}})$$

$$\geq\lim_{m\to\infty}(b\frac{m}{\log m}+\frac{r_m-1}{\log m})(b+\frac{r_m}{m})>b^2>0.$$

Hence

$$\lim_{m\to\infty}\varepsilon_m(n)\left[\{1+\frac{\log\log m-1}{\log m}+O\left(\frac{\log\log m}{(\log m)^2}\right)\}\right](1+\frac{0.0066788}{\log p_m})=0>b^2>0.$$

It is contradictory. Thus, we prove that $\lim_{m\to\infty}\frac{K_0}{m}=0$.

**Proof of inference 2** since $1=\alpha_1=\alpha_2=\cdots=\alpha_m$, $\varepsilon_m(n)=0$, $\lim_{n\to\infty}\varepsilon_m(n)=0$, $A=K_0=1$. According to Theorem 4, the conclusion is correct. □

# 4 Proof of lemmas

## 4.1 Proof of Lemma 2.1

Simple calculation shows

$$\frac{\partial f}{\partial q_k} = e^{\gamma} \frac{\alpha_k q_k^{\alpha_k-1} \prod_{h \neq k, 1 \leq h \leq m} q_h^{\alpha_h}}{\left(\log \prod_{h=1}^{m} q_h^{\alpha_h}\right) \prod_{h=1}^{m} q_h^{\alpha_h}} - \prod_{h \neq k, 1 \leq h \leq m} \frac{q_h - \frac{1}{q_h^{\alpha_h}}}{q_h - 1} \left[\frac{1 + \frac{\alpha_k}{q_k^{\alpha_k+1}}}{q_k - 1} - \frac{q_k - \frac{1}{q_k^{\alpha_k}}}{(q_k - 1)^2}\right]$$

$$= e^{\gamma} \frac{\alpha_k}{q_k \left(\log \prod_{h=1}^{m} q_h^{\alpha_h}\right)}$$

$$- \prod_{h \neq k, 1 \leq h \leq m} \frac{q_h - \frac{1}{q_h^{\alpha_h}}}{q_h - 1} \left(\frac{\alpha_k}{q_k^{\alpha_k}} + \frac{1}{q_k^{\alpha_k}} - 1 - \frac{\alpha_k}{q_k^{\alpha_k+1}}\right) \frac{1}{(q_k - 1)^2}.$$

In the following we firstly prove

$$\frac{\alpha_k}{q_k^{\alpha_k}} + \frac{1}{q_k^{\alpha_k}} - 1 - \frac{\alpha_k}{q_k^{\alpha_k+1}} < 0,$$

namely

$$(\alpha_k + 1)q_k - q_k^{\alpha_k+1} - \alpha_k < 0.$$

Let

$$g(q_k) = (\alpha_k + 1)q_k - q_k^{\alpha_k+1} - \alpha_k.$$

Then

$$g'(q_k) = (\alpha_k + 1)(1 - q_k^{\alpha_k}).$$

Since $\alpha_k \geq 1$ and $q_k \geq 2$, we have $g'(q_k) < 0$ and $g(q_k)$ is a monotone decreasing function of $q_k$. Therefore

$$g(q_k) \leq g(2) = \alpha_k + 2 - 2^{\alpha_k+1}.$$

Let

$$t(\alpha_k) = \alpha_k + 2 - 2^{\alpha_k+1}.$$

Then $t'(\alpha_k) = 1 - 2^{\alpha_k+1}\log 2$. By $\alpha_k \geq 1$, we have $t'(\alpha_k) < 0$, and $t(\alpha_k) \leq t(1) = 1 + 2 - 2^{1+1} < 0$. In summary, $g(q_k) < 0$. Hence $\dfrac{\partial f}{\partial q_k} > 0 (k = 1, 2, \cdots, m)$. Thus $f(q_1, q_2, \cdots, q_m)$ has no extreme points in

$$D: \{(q_1, q_2, \cdots, q_m) | 2 \leq q_1 < q_2 < \cdots < q_m,\ q_h \in P,\ h = 1, 2, \cdots, m\},$$

and the minimum of $f(q_1, q_2, \cdots, q_m)$ can be found at the boundary $\partial D$. Since $\dfrac{\partial f}{\partial q_k} > 0$ $(k = 1, 2, \cdots, m)$. Therefore $f(q_1, \cdots, q_m)$ is monotone-increasing with respect to any $q_k$. Furthermore $q_1 < q_2 < \cdots < q_m$, thus $f(q_1, q_2, \cdots, q_m)$ has its minimum at the boundary point

$$(q_1, q_2, \cdots, q_m) = (2, 3, 5, \cdots, p_m),$$

namely

$$\min_{(q_1,q_2,\cdots,q_m)\in D} f(q_1, q_2, \cdots, q_m) = f(2, 3, 5, \cdots, p_m) = f(p_1, p_2, \cdots, p_m). \quad \square$$

**4.2 Proof of lemma 2.3**

If $\lim\limits_{m \to \infty} \varepsilon_m(n) \neq 0$, then $\lim\limits_{m \to \infty} \varepsilon_m(n) = +\infty$ or $\lim\limits_{m \to \infty} \varepsilon_m(n) \neq +\infty$. When $\lim\limits_{m \to \infty} \varepsilon_m(n) = +\infty$, obviously there exists constant $a > 0$, such that $\lim\limits_{m \to \infty} \varepsilon_m(n) > 2a > 0$. Let

$$\varepsilon_s = \varepsilon_s(n) = \left(\log \prod_{k=1}^{s} p_k^{\alpha_k}\right)\left(\log \prod_{k=1}^{s} p_k\right)^{-1} - 1 (s = 1, 2, \cdots, m).$$

In the case of $\lim\limits_{m \to \infty} \varepsilon_m(n) \neq +\infty$, we prove $\varepsilon_1 \geq \varepsilon_2 \geq \cdots \geq \varepsilon_k \geq \cdots \geq \varepsilon_m \geq 0$.

Since $\alpha_1 \geq \alpha_2 \geq \cdots \geq \alpha_m \geq 1$, $\alpha_1, \alpha_2, \cdots, \alpha_m \in \mathbb{N}$, for $k: 1 \leq k < m$,

$$\prod_{i=1}^{k} p_i^{\alpha_i} = \left(\prod_{i=1}^{k} p_i\right)^{1+\varepsilon_k}, \quad \prod_{i=1}^{k+1} p_i^{\alpha_i} = \left(\prod_{i=1}^{k+1} p_i\right)^{1+\varepsilon_{k+1}}.$$

Then

$$\left(\prod_{i=1}^{k} p_i\right)^{1+\varepsilon_k} p_{k+1}^{\alpha_{k+1}} = \left(\prod_{i=1}^{k} p_i^{\alpha_i}\right) p_{k+1}^{\alpha_{k+1}} = \prod_{i=1}^{k+1} p_i^{\alpha_i}$$

$$= \left(\prod_{i=1}^{k+1} p_i\right)^{1+\varepsilon_{k+1}} = \left(\prod_{i=1}^{k} p_i\right)^{1+\varepsilon_{k+1}} p_{k+1}^{1+\varepsilon_{k+1}},$$

namely

$$\left(\prod_{i=1}^{k} p_i\right)^{1+\varepsilon_{k+1}} p_{k+1}^{1+\varepsilon_{k+1}} = \left(\prod_{i=1}^{k} p_i\right)^{1+\varepsilon_k} p_{k+1}^{\alpha_{k+1}}. \tag{4.1}$$

On the other hand

$$\alpha_1 \geqslant \alpha_2 \geqslant \cdots \geqslant \alpha_k \geqslant \alpha_{k+1}.$$

Hence

$$\left(\prod_{i=1}^{k+1} p_i\right)^{1+\varepsilon_{k+1}} = \prod_{i=1}^{k+1} p_i^{\alpha_i} \geqslant \prod_{i=1}^{k+1} p_i^{\alpha_{k+1}} = \left(\prod_{i=1}^{k+1} p_i\right)^{\alpha_{k+1}}.$$

Thus $1+\varepsilon_{k+1} \geqslant \alpha_{k+1}$, according to (4.1),

$$\left(\prod_{i=1}^{k} p_i\right)^{1+\varepsilon_{k+1}} p_{k+1}^{1+\varepsilon_{k+1}} = \left(\prod_{i=1}^{k} p_i\right)^{1+\varepsilon_k} p_{k+1}^{\alpha_{k+1}} \leqslant \left(\prod_{i=1}^{k} p_i\right)^{1+\varepsilon_k} p_{k+1}^{1+\varepsilon_{k+1}}.$$

Therefore, $\varepsilon_k \geqslant \varepsilon_{k+1}$, thus $\varepsilon_1 \geqslant \varepsilon_2 \geqslant \cdots \geqslant \varepsilon_k \geqslant \cdots \geqslant \varepsilon_m$.

Since

$$\alpha_1 \geqslant \alpha_2 \geqslant \cdots \geqslant \alpha_m \geqslant 1, \quad \left(\prod_{i=1}^{k} p_i\right)^{1+\varepsilon_k} = \prod_{i=1}^{k} p_i^{\alpha_i}, \quad 1 \leqslant k \leqslant m.$$

Then $\varepsilon_1 \geqslant \varepsilon_2 \geqslant \cdots \geqslant \varepsilon_k \geqslant \cdots \geqslant \varepsilon_m \geqslant 0$, by $\lim_{m\to\infty} \varepsilon_m \neq 0$, according to the Cauchy convergence criterion, there exists constant $a > 0$, such that $\lim_{m\to\infty} \varepsilon_m > 2a > 0$. Lemma 2.5 is proved. □

**Note 3** for any given $n = \prod_{i=1}^{m} p_i^{\alpha_i} = \left(\prod_{i=1}^{m} p_i\right)^{1+\varepsilon_m(n)}$ ($\alpha_1 \geqslant \alpha_2 \geqslant \cdots \geqslant \alpha_m \geqslant 1$), we always have $\varepsilon_1 \geqslant \varepsilon_2 \geqslant \cdots \geqslant \varepsilon_k \geqslant \cdots \geqslant \varepsilon_m$ and $\varepsilon_m(n) \in [0, +\infty)$. For example, let $n = \prod_{k=1}^{m} p_k^{(m-k+1)!}$, then

$$\varepsilon_1 = \frac{m!\log 2}{\log 2} - 1 > \varepsilon_2 = \frac{m!\log 2 + (m-1)!\log 3}{\log 2 + \log 3} - 1$$

$$> \varepsilon_3 = \frac{m!\log 2 + (m-1)!\log 3 + (m-2)!\log 5}{\log 2 + \log 3 + \log 5} - 1$$

$$> \varepsilon_4 > \cdots > \varepsilon_m = \frac{\sum_{k=1}^{m}(m+1-k)!\log p_k}{\sum_{k=1}^{m}\log p_k} - 1 > \frac{m!\log p_1}{m\log p_m} \to +\infty (m \to \infty).$$

When $n = \prod_{k=1}^{m} p_k = \prod_{k=1}^{m} p_k^{1+\varepsilon_m(n)}$, then $\varepsilon_1 = \varepsilon_2 = \cdots = \varepsilon_m = 0$.

### 4.3 Proof of Lemma 2.4

**Proposition 4.1** If $0 < \varepsilon < \frac{1}{5}$ and $c$ is a positive constant, then

$$\lim_{x \to +\infty} (\log x) O\left(e^{-c\log^{\frac{3}{5}-\varepsilon} x}\right) = 0,$$

where $O\left(e^{-c\log^{\frac{3}{5}-\varepsilon} x}\right)$ is defined as: there exist positive constants $A$ and $B$, for any $x > B$, such that $\left|O\left(e^{-c\log^{\frac{3}{5}-\varepsilon} x}\right)\right| \leqslant A e^{-c\log^{\frac{3}{5}-\varepsilon} x}$.

**Proof** According to the definition of $O\left(e^{-c\log^{\frac{3}{5}-\varepsilon} x}\right)$, there exists positive constants $A$, $B$ and for all $x > B$, such that

$$-A e^{-c\log^{\frac{3}{5}-\varepsilon} x} \leqslant O\left(e^{-c\log^{\frac{3}{5}-\varepsilon} x}\right) \leqslant A e^{-c\log^{\frac{3}{5}-\varepsilon} x}.$$

According to L' Hospital's rule

$$\lim_{x \to +\infty} (\log x) A e^{-c\log^{\frac{3}{5}-\varepsilon} x} = \lim_{x \to +\infty} \frac{A \log x}{e^{c\log^{\frac{3}{5}-\varepsilon} x}} = \lim_{x \to +\infty} \frac{A \log^{\frac{2}{5}+\varepsilon} x}{c(\frac{3}{5}-\varepsilon) e^{c\log^{\frac{3}{5}-\varepsilon} x}}.$$

Let $\mu = \log x$:

x

$$\lim_{x \to +\infty} \frac{A \log^{\frac{2}{5}+\varepsilon} x}{c(\frac{3}{5}-\varepsilon)e^{c\log^{\frac{3}{5}-\varepsilon} x}} = \lim_{\mu \to +\infty} \frac{A\mu^{\frac{2}{5}+\varepsilon}}{c(\frac{3}{5}-\varepsilon)e^{c\mu^{\frac{3}{5}-\varepsilon}}}.$$

Applying L' Hospital's rule gets

$$\lim_{x \to +\infty}(\log x)Ae^{-c\log^{\frac{3}{5}-\varepsilon} x} = \lim_{\mu \to +\infty} \frac{A\mu^{\frac{2}{5}+\varepsilon}}{c(\frac{3}{5}-\varepsilon)e^{c\mu^{\frac{3}{5}-\varepsilon}}} = \lim_{\mu \to +\infty} \frac{A(\frac{2}{5}+\varepsilon)\mu^{2\varepsilon-\frac{1}{5}}}{c^2(\frac{3}{5}-\varepsilon)^2 e^{c\mu^{\frac{3}{5}-\varepsilon}}}.$$

If $2\varepsilon - \frac{1}{5} \leqslant 0$, then $\lim_{\mu \to +\infty} \dfrac{A(\frac{2}{5}+\varepsilon)\mu^{2\varepsilon-\frac{1}{5}}}{c^2(\frac{3}{5}-\varepsilon)^2 e^{c\mu^{\frac{3}{5}-\varepsilon}}} = 0$.

If $2\varepsilon - \frac{1}{5} > 0$, since $0 < \varepsilon < \frac{1}{5}$, $3\varepsilon - \frac{4}{5} < 0$, applying L' Hospital's rule again gets

$$\lim_{\mu \to +\infty} \frac{A(\frac{2}{5}+\varepsilon)\mu^{2\varepsilon-\frac{1}{5}}}{c^2(\frac{3}{5}-\varepsilon)^2 e^{c\mu^{\frac{3}{5}-\varepsilon}}} = \lim_{\mu \to +\infty} \frac{A(\frac{2}{5}+\varepsilon)(2\varepsilon-\frac{1}{5})\mu^{3\varepsilon-\frac{4}{5}}}{c^3(\frac{3}{5}-\varepsilon)^3 e^{c\mu^{\frac{3}{5}-\varepsilon}}} = 0.$$

Similarly $\lim_{x \to +\infty}(\log x)(-Ae^{-c\log^{\frac{3}{5}-\varepsilon} x}) = 0$. According to the squeeze theorem

$$\lim_{x \to +\infty}(\log x)O\left(e^{-c\log^{\frac{3}{5}-\varepsilon} x}\right) = 0.$$

The proposition is proved. □

The following result is well known to us:

**Proposition 4.2([10])**

$$\prod_{k=1}^{m}\left(1 - \frac{1}{p_k}\right) = \frac{1}{e^{\gamma}\log x}\left\{1 + O\left(e^{-c\log^{\frac{3}{5}-\varepsilon} x}\right)\right\} \quad (x \geqslant 2).$$

Here $\gamma$ is Euler's constant, $c$ is a constant greater than $0$, and $\varepsilon$ is a positive number that can be arbitrarily small.

To prove Lemma 2.4, the following proposition is required.

**Proposition 4.3** For any $m \in \mathbb{N}$, $\alpha_k \geqslant 1(k = 1, 2, \cdots, m)$, such that

$$\prod_{k=1}^{m}\frac{p_k-\frac{1}{p_k^{\alpha_k}}}{p_k-1}\leqslant e^{\gamma}\log p_m\left\{1+O\left(e^{-c\log^{\frac{3}{5}-\varepsilon}p_m}\right)\right\}^{-1}, \qquad (4.2)$$

where $c$ is a constant greater than $0$ and $\varepsilon$ is a positive that can be arbitrarily small.

**Proof** Since $m\in\mathbb{N}$, $\alpha_k\geqslant 1(k=1,2,\cdots,m)$,

$$\prod_{k=1}^{m}\frac{p_k-\frac{1}{p_k^{\alpha_k}}}{p_k-1}<\prod_{k=1}^{m}\frac{p_k}{p_k-1}=\left\{\prod_{k=1}^{m}\left(1-\frac{1}{p_k}\right)\right\}^{-1}.$$

According to proposition 4.2,

$$\left\{\prod_{k=1}^{m}\left(1-\frac{1}{p_k}\right)\right\}^{-1}=\left(\frac{1}{e^{\gamma}\log p_m}\left\{1+O\left(e^{-c\log^{\frac{3}{5}-\varepsilon}p_m}\right)\right\}\right)^{-1}$$

$$=\left(e^{\gamma}\log p_m\right)\left\{1+O\left(e^{-c\log^{\frac{3}{5}-\varepsilon}p_m}\right)\right\}^{-1}.$$

Thus, (4.2) is true. □

**Proposition 4.4[11]**   For any $\alpha_k\geqslant 1$,

$$\log\prod_{k=1}^{m}p_k^{\alpha_k}=\sum_{k=1}^{m}\alpha_k\log p_k\geqslant\sum_{p\leq p_m}\log p$$

$$=\vartheta(p_m)=p_m+O\left(\frac{p_m}{\log p_m}\right).$$

**Proof** From [8, 12],

$$\frac{x}{\log x}(1+\frac{1}{2\log x})<\pi(x)<\frac{x}{\log x}(1+\frac{3}{2\log x}),$$

where $\pi(x)$ be the number of primes not exceeding $x$. Hence

$$\pi(x)\log x=x+O\left(\frac{x}{\log x}\right).$$

On the other hand

$$\pi(x)=\frac{\vartheta(x)}{\log x}+O\left(1+\frac{x}{(\log x)^2}\right).$$

Thus

$$\vartheta(x) = x + O\left(\frac{x}{\log x}\right).$$

The proposition is proved for $x = p_m$. □

**Proof of lemma 2.4** According to (4.2),
$$\prod_{k=1}^{m} \frac{p_k - \frac{1}{p_k^{\alpha_k}}}{p_k - 1} \leqslant \left(e^{\gamma} \log p_m\right)\left\{1 + O\left(e^{-c\log^{\frac{3}{5}-\varepsilon} p_m}\right)\right\}^{-1}.$$

According to Proposition 4.4,
$$\vartheta(p_m) = p_m + O\left(\frac{p_m}{\log p_m}\right).$$

$$\lim_{m \to +\infty}\left(e^{\gamma} \log\log \prod_{k=1}^{m} p_k - \prod_{k=1}^{m} \frac{p_k - \frac{1}{p_k^{\alpha_k}}}{p_k - 1}\right)$$

$$= \lim_{m \to +\infty}\left(e^{\gamma} \log(\vartheta(p_m)) - \prod_{k=1}^{m} \frac{p_k - \frac{1}{p_k^{\alpha_k}}}{p_k - 1}\right)$$

$$\geqslant \lim_{m \to +\infty}\left(e^{\gamma} \log(\vartheta(p_m))\right)$$

$$-\left(e^{\gamma} \log p_m\right)\left\{1 + O\left(e^{-c\log^{\frac{3}{5}-\varepsilon} p_m}\right)\right\}^{-1}\right)$$

$$= \lim_{m \to +\infty}\left(e^{\gamma} \log\left(p_m + O\left(\frac{p_m}{\log p_m}\right)\right)\right.$$

$$\left.-\left(e^{\gamma} \log p_m\right)\left\{1 + O\left(e^{-c\log^{\frac{3}{5}-\varepsilon} p_m}\right)\right\}^{-1}\right)$$

$$= e^{\gamma} \lim_{p_m \to +\infty}\left(\log\left(p_m + O\left(\frac{p_m}{\log p_m}\right)\right)\right.$$

$$\left.-(\log p_m)\left\{1 + O\left(e^{-c\log^{\frac{3}{5}-\varepsilon} p_m}\right)\right\}^{-1}\right)$$

$$= e^{\gamma} \lim_{p_m \to +\infty}\left(\log p_m + \log\left(1 + \frac{O\left(\frac{p_m}{\log p_m}\right)}{p_m}\right)\right.$$

$$\left.-(\log p_m)\left\{1 + O\left(e^{-c\log^{\frac{3}{5}-\varepsilon} p_m}\right)\right\}^{-1}\right)$$

$$= e^\gamma \lim_{p_m \to +\infty} \left( \frac{(\log p_m) O\left(e^{-c \log^{\frac{3}{5}-\varepsilon} p_m}\right)}{1 + O\left(e^{-c \log^{\frac{3}{5}-\varepsilon} p_m}\right)} + \log\left(1 + \frac{O\left(\frac{p_m}{\log p_m}\right)}{p_m}\right) \right).$$

Since $\varepsilon$ is an arbitrarily small number, according to Proposition 4.1,

$$e^\gamma \lim_{p_m \to +\infty} \left( \frac{(\log p_m) O\left(e^{-c \log^{\frac{3}{5}-\varepsilon} p_m}\right)}{1 + O\left(e^{-c \log^{\frac{3}{5}-\varepsilon} p_m}\right)} + \log\left(1 + \frac{O\left(\frac{p_m}{\log p_m}\right)}{p_m}\right) \right)$$

$$= \frac{e^\gamma \lim_{p_m \to +\infty} (\log p_m) O\left(e^{-c \log^{\frac{3}{5}-\varepsilon} p_m}\right)}{\lim_{p_m \to +\infty}\left(1 + O\left(e^{-c \log^{\frac{3}{5}-\varepsilon} p_m}\right)\right)}$$

$$+ e^\gamma \lim_{p_m \to +\infty} \log\left(1 + \frac{O\left(\frac{p_m}{\log p_m}\right)}{p_m}\right) = e^\gamma \left(\frac{0}{1} + 0\right) = 0.$$

(2.3) holds, and the Lemma is proved. □

**Acknowledgments**

The author thanks Dr. Ru Zhu at Graceland University, Dr. Lei Wang, Prof. Yuwu Yao at Hefei University, Dr. Hongbin Fang at Georgetown University and Dr. Wenjin Lv at Xi'an Jiaotong-liverpool University for helpful discussions and translation.

**References**

[1] G. F. B. Riemann. Ueber die Anzahl der Primzalen unter einer gegebenen Grösse, Monatsber. Akad. Berlin (1859): 671-680.

[2] G. Robin. Grandes valeurs de la function somme des diviseurs et hypothèse de Riemann, J. Math. Pures Appl. **63**(1984): 184-213.

[3] S. Nazardonyavi, S. Yakubovich. Superabundant numbers, their subsequences and the Riemann hypothesis. arXiv:1211.2147v3 [math.NT] 26 Feb 2013.

[4] Y. Zhu. On A Riemann Hypothesis Related Inequality. Journal of Hefei University, Vol.26 No.1 (2016,1): 1-8.

[5] Y. ZHU Yuyang, Study on a Riemann Hypothesis related Inequality. http://www.paper.edu.cn/releasepaper/content/201512-541.

[6] P. Solé and Y. Zhu. An Asymptotic Robin Inequality. INTEGERS, #A81, 16 (2016).

[7] P. Dusart. Intégalités explicites pour $\psi(x), \theta(x), \pi(x)$ et les nombres premiers. C. R. Math. Rep. Acad. Sci. Canada 21(1999), 53-59.

[8] L. Schoenfeld. Sharper bounds for the Chebyshev functions $\vartheta(x)$ and $\psi(x)$, II, Math. Comp. 1976, 30, 337-360.

[9] P. Dusart. The $k^{th}$ prime is greater than $k(lnk+lnlnk-1)$ for $k \geq 2$. Math. Comp. 68

(1999), 411-415.

[10] A. A. Kapauyóa. Basic analytic number theory, HAYKA, 1975, 65-75. (in Russia)

[11] G. Tenenbaum. Introduction to analytic and probability number theory, Cambridge University Press. 1998.

[12] J. B. Rosser, L. Schoenfeld. Approximate formulas for some functions of prime numbers, Illinois J. Math. 1962, 6, 64-94.